\documentclass[aop,seceqn,citesort,MSNbibl,dvips]{arximspdf}
\usepackage{mathrsfs}

\doi{10.1214/10-AOP600}
\volume{39}
\issue{4}
\pubyear{2011}
\firstpage{1449}
\lastpage{1467}

\makeatletter

\newtheorem{theorem}{Theorem}[section]
\newtheorem{cor}[theorem]{Corollary}
\newtheorem{lem}[theorem]{Lemma}

\newproclaim{Remark}{Remark}[section]

\newcommand{\R}{\mathbb R}
\newcommand{\ff}{\frac}
\newcommand{\B}{\mathscr{B}}
\newcommand{\kk}{\kappa}
\newcommand{\DD}{\Delta}
\newcommand{\rr}{\rho}

\newcommand{\GG}{\Gamma}
\newcommand{\nn}{\nabla}
\newcommand{\pp}{\partial}
\newcommand{\ddd}{{{d}}}
\newcommand{\bb}{\beta}

\newcommand{\si}{\sigma}

\newcommand{\F}{\mathscr{F}}

\newcommand{\OO}{\Omega}

\newcommand{\Ric}{\operatorname{Ric}}

\newcommand{\II}{\mathbb I}

\newcommand{\E}{\mathbb E}

\newcommand{\Q}{\mathbb Q}

\makeatother

\begin{document}
\begin{frontmatter}

\title{Harnack inequality for SDE with multiplicative noise and
extension to Neumann semigroup on~nonconvex
manifolds\thanksref{T1}}
\runtitle{Harnack inequalities and applications}

\thankstext{T1}{Supported in part by WIMCS, NNSFC (10721091) and the
973-Project.}

\begin{aug}
\author[A]{\fnms{Feng-Yu} \snm{Wang}\corref{}\ead[label=e1]{wangfy@bnu.edu.cn}\ead[label=e2]{F.Y.Wang@swansea.ac.uk}}
\runauthor{F.-Y. Wang}
\affiliation{Beijing Normal University and Swansea University}
\address[A]{School of Mathematical Sciences\\
Laboratory of Mathematics\\
\quad and Complex Systems\\
Beijing Normal University\\
Beijing 100875\\
China\\
\printead{e1}\\
and\\
Department of Mathematics\\
Swansea University\\
Singleton Park, SA2 8PP\\
United Kingdom\\
\printead{e2}} 
\end{aug}

\received{\smonth{11} \syear{2009}}
\revised{\smonth{9} \syear{2010}}

%
\begin{abstract}
By constructing a coupling with unbounded time-dependent drift,
dimension-free Harnack inequalities are established for a large class
of stochastic differential equations with multiplicative noise. These
inequalities are applied to the study of heat kernel upper bound and
contractivity properties of the semigroup. The main results are also
extended to reflecting diffusion processes on Riemannian manifolds with
nonconvex boundary.
\end{abstract}

%

\end{frontmatter}

\section{Introduction}\label{sec1}

Consider the following SDE on $\R^d$:
%
%
\begin{equation}\label{1.1}
\ddd X_t = \si(t,X_t)\,\ddd B_t + b(t, X_t)\,\ddd t,
\end{equation}
where $B_t$ is the $d$-dimensional Brownian motion
on a complete filtered probability space $(\OO, \{\F_t\}_{t\ge0},
\mathbb P)$, and
\[
\si\dvtx[0,\infty)\times\R^d\times\OO\to\R^d\otimes\R^d,\qquad
b\dvtx[0,\infty)\times\R^d\times\OO\to\R^d
\]
are progressively
measurable and continuous in the second variable. Throughout the
paper, we assume that for any $X_0\in\R^d$ the equation
(\ref{1.1}) has a unique strong solution which is nonexplosive and
continuous in $t$.

Let $X_t^x$ be the solution to (\ref{1.1}) for $X_0=x$. We aim to
establish the Harnack inequality for the operator $P_t$:
\[
P_t f(x):= \E f(X_t^x),\qquad t\ge0, x\in\R^d, f\in
\B_b^+(\R^d),
\]
where $\B_b^+(\R^d)$ is the class of all bounded
nonnegative measurable functions on~$\R^d$.
To this end, we shall make use of
the following assumptions.

\begin{longlist}[(A2)]
\item[(A1)] There exists an increasing function $K\dvtx[0,\infty)\to
\R$
such that almost surely
\begin{eqnarray*}
&&\|\si(t,x)-\si(t,y)\|_{\mathrm{HS}}^2 + 2 \langle b(t,x)-b(t,y),
x-y\rangle\\
&&\qquad\le
K_t|x-y|^2,\qquad x,y\in\R^d, t\ge0.
\end{eqnarray*}
\item[(A2)] There exists a
decreasing function $\lambda\dvtx[0,\infty)\to(0,\infty)$ such that
almost surely
\[
\si(t,x)^*\si(t,x)\ge\lambda_t^2 I,\qquad x\in\R^d, t\ge0.
\]
\item[(A3)] There exists an increasing function $\delta\dvtx[0,\infty
)\to
(0,\infty)$ such that almost surely
\[
\bigl|\bigl(\si(t, x)- \si(t,y)\bigr)(x-y)\bigr|\le\delta_t |x-y|,\qquad x,y\in\R^d, t\ge
0.
\]
\item[(A4)] For $n\ge1$, there exists a constant $c_n>0$ such
that almost surely
\[
\|\si(t,x)-\si(t,y)\|_{\mathrm{HS}} + |b(t,x)-b(t,y)|\le c_n |x-y|,\qquad |x|,
|y|, t\le n.
\]
\end{longlist}

It is well known that (A1) ensures the uniqueness of the
solution to (\ref{1.1}) while (A4) implies the existence and the
uniqueness of the strong solution (see, e.g.,~\cite{FZ} and references within
for weaker conditions). On the other hand, if $b$ and $\si$ depend only
on the
variable $x\in\R^d$, then their continuity in $x$ implies the
existence of weak solutions (see~\cite{IW}, Theorem 2.3), so that by the
Yamada--Watanabe principle~\cite{YW}, the uniqueness ensured by (A1)
implies the existence and
uniqueness of the strong solution.

Note that if $\si(t,x)$ and $b(t,x)$ are deterministic and
independent of $t$, then the solution is a time-homogeneous Markov
process generated by
\[
L:= \ff1 2 \sum_{i,j=1}^d a_{ij} \,\pp_i\,\pp_j +\sum_{i=1}^d b_i\,
\pp_i,
\]
where $a:=\si\si^*$. If further more $\si$ and $b$ are
smooth, we may consider the Bakry--Emery curvature condition
\cite{BE}:
%
%
\begin{equation}\label{CC}\GG_2(f,f)\ge-K\GG(f,f),\qquad f\in
C^\infty(\R^d),
\end{equation}
for some constant $K\in\R$, where
\begin{eqnarray*}
\GG(f,g)&:= &\ff1 2 \sum_{i,j=1}^d a_{ij}
(\pp_i f)(\pp_j g),\qquad f,g\in C^1(\R^d),\\
\GG_2(f,f)&:=& \tfrac1 2 L\GG(f,f)- \GG(f, Lf),\qquad f\in
C^\infty(\R^d).
\end{eqnarray*}
According to~\cite{W97}, Lemma 2.2, and~\cite{W04}, Theorem 1.2, the
curvature condition
(\ref{CC}) is equivalent to the dimension-free Harnack inequality
\begin{eqnarray}
&&(P_t f(x))^p \le(P_t f^p(y))
\exp\biggl[\ff{p\rr_a(x,y)^2}{2(p-1) (1-{e}^{-2Kt})}\biggr],\nonumber\\
&&\eqntext{t\ge0,
p>1, f\in\B_b^+(\R^d), x,y\in\R^d,}
\end{eqnarray}
where
\[
\rr_a(x,y):= \sup\{|f(x)-f(y)|\dvtx f\in C^1(\R^d), \GG(f,f)\le1\}
,\qquad x,y\in\R^d.
\]
This type of inequality has
been extended and applied to the study of heat kernel (or transition
probability) and contractivity properties for diffusion semigroups,
see
\mbox{\cite{A98,RW03a,ATW}} for diffusions on manifolds with possibly
unbounded below curvature,~\cite{W07,LW08} for stochastic
generalized porous media and fast diffusion equations, and
\cite{AK,AZ,DRW,OY,Liu,RW03,ERS,Zh} for the study of some other SPDEs with
additive noise.

If $\si$ depends on $x$, however, it is normally very hard to
verify the curvature condition (\ref{CC}), which
depends on second order derivatives of $a^{-1}$, the inverse matrix
of $a$. This is the main reason why existing results on the
dimension-free Harnack inequality for SPDEs are only proved for
the additive noise case.

In this paper, we shall use the coupling argument developed in~\cite{ATW},
which will allow us to establish Harnack inequalities for $\si(t,x)$
depending on $x$.
This method has also been applied to the study of SPDEs in the above
mentioned references. To see the difficulty in the study for
$\si(t,x)$ depending on $x$, let us briefly recall the main idea of this
argument.

To explain the main idea of the coupling, we first consider the easy
case where $\si$ and $b$ are
independent of the second variable. For $x\ne y$ and $T>0$, let
$X_t$ solve (\ref{1.1}) with $X_0=x$ and $Y_t$ solve
\[
\ddd Y_t = \si(t) \,\ddd B_t +b(t)\,\ddd t
+\ff{|x-y|(X_t-Y_t)}{T|X_t-Y_t|}\,\ddd t,\qquad Y_0=y.
\]
Then $Y_t$ is
well defined up to the coupling time
\[
\tau:=\inf\{t\ge0\dvtx X_t=Y_t\}.
\]
Let $X_t=Y_t$ for $t\ge\tau$.
We have
\[
\ddd|X_t-Y_t|= -\ff{|x-y|}T\,\ddd t,\qquad t\le\tau.
\]
This implies
$\tau=T$ and hence, $X_T=Y_T$. On the other hand, by the Girsanov
theorem we have
\[
P_T f(y)= \E[ R f(Y_T)]
\]
for
\begin{eqnarray*}
R&:=& \exp\biggl[-\ff{|x-y|}T \int_0^T \ff{\langle\si(t)^{-1} (X_t-Y_t),
\ddd
B_t\rangle}{|X_t-Y_t|} \\
&&\hspace*{20.2pt}{}- \ff{|x-y|^2}{2 T^2} \int_0^T \ff{|\si(t)^{-1}
(X_t-Y_t)|^2}{|X_t-Y_t|^2}\,\ddd t\biggr].
\end{eqnarray*}
Therefore,
\[
(P_T f(y))^p = (\E[R f(X_T)])^p \le(P_T f^p(x)) \bigl(\E
R^{p/(p-1)}\bigr)^{p-1}.
\]
Since by (A1) and (A2) it is easy to
estimate moments of $R$, the desired Harnack inequality follows
immediately.

In general, if $\si(t,x)$ depends on $x$, then the process
$X_t-Y_t$ contains a nontrivial martingale term, which cannot be
dominated by and bounded drift. So, in this case, any additional
bounded drift put in the equation for $Y_t$ is not enough to make
the coupling successful before a fixed time $T$. This is the main
difficulty to establish the Harnack inequality for diffusion
semigroups with nonconstant diffusion coefficient.

In this paper, under assumptions (A1) and (A2), we are able to
constructed a coupling with a drift which is unbounded around a
fixed time $T$, such that the coupling is successful before $T$.
In this case, the corresponding exponential martingale has finite
entropy such that the log-Harnack inequality holds; if further more
(A3) holds then the exponential martingale is $L^p$-integrable
for some $p>1$ such that the Harnack inequality with power holds.
More precisely, we have the following result.
\begin{theorem}\label{T1.1} Let $\si(t,x)$ and $b(t,x)$ either be
deterministic and independent of $t$, or satisfy \textup{(A4)}.
\begin{longlist}[(2)]
\item[(1)] If \textup{(A1)} and \textup{(A2)} hold, then
\begin{eqnarray}
&&P_T\log f(y)\le\log P_T f(x) + \ff{ K_T|x-y|^2}{2\lambda_T^2 (1-{e}
^{-K_TT})},\nonumber\\
&&\eqntext{f\ge1, x,y\in\R^d, T>0.}
\end{eqnarray}
\item[(2)] If
\textup{(A1)}, \textup{(A2)} and \textup{(A3)} hold, then for $p>(1+\ff
{\delta_T}{\lambda_T})^2$
and $\delta_{p,T}:=\max\{\delta_T$, $\ff{\lambda_T} 2 (\sqrt p-1)\}
$, the Harnack
inequality
\[
(P_T f(y))^p \le(P_T f^p(x))\exp\biggl[\ff{K_T\sqrt p (\sqrt p-1)
|x-y|^2}{4\delta_{p,T} [(\sqrt p-1)\lambda_T -\delta_{p,T}]
(1-{e}^{-K_TT})}\biggr]
\]
holds for all $T>0, x,y\in\R^d$ and $f\in
\B_b^+(\R^d)$.
\end{longlist}
\end{theorem}

Theorem~\ref{T1.1}(1) generalizes a recent
result in~\cite{RW09} on the log-Harnack inequality by using the
gradient estimate on $P_t$.

Let $p_t(x,y)$ be the density of $P_t$ w.r.t. a Radon measure $\mu$.
Then according to~\cite{WN}, Proposition 2.4, the above log-Harnack
inequality and
Harnack inequality are equivalent to the following heat kernel
inequalities, respectively:
%
%
\begin{eqnarray}\label{Heat1}
&&\int_{\R^d}
p_T(x,z)\log\ff{p_T(x,z)}{p_T(y,z)}\mu(\ddd z)
\le\ff{
K|x-y|^2}{2\lambda_T^2 (1-{e}^{-K_TT})},\nonumber\\[-8pt]\\[-8pt]
&&\eqntext{x,y\in\R^d,
T>0,}
\end{eqnarray}
and
%
%
\begin{eqnarray}\label{Heat2}
&&\int_{\R^d} p_T (x,z)\biggl(\ff{p_t(x,z)}{p_t(y,z)}
\biggr)^{1/(p-1)}\mu(\ddd z)\nonumber\\
&&\qquad\le\exp\biggl[\ff{K_T\sqrt p |x-y|^2}{4\delta_{p,T} (\sqrt p+1)[(\sqrt
p-1)\lambda_T -\delta_{p,T}] (1-{e}^{-K_TT})}\biggr],\\
&&\eqntext{x,y\in\R^d,
T>0.}
\end{eqnarray}
So, the following is a direct
consequence of Theorem~\ref{T1.1}.
\begin{cor} \label{C1.0} Let $\si(t,x)$ and $b(t,x)$ either be
deterministic and independent of $t$, or satisfy \textup{(A4)}. Let $P_t$
have a strictly positive density $p_t(x,y)$ w.r.t. a Radon measure $\mu
$. Then
\textup{(A1)} and \textup{(A2)} imply (\ref{Heat1}), while \textup
{(A1)}--\textup{(A3)} imply
(\ref{Heat2}).
\end{cor}

Next, by standard applications of the Harnack inequality with power, we
have the
following consequence of Theorem~\ref{T1.1} on contractivity properties
of $P_t$.
\begin{cor}\label{C1.2} Let $\si(t,x)$ and $b(t,x)$ be deterministic
and independent of $t$, such that \textup{(A1)--(A3)} hold for constant
$K,\lambda$ and $\delta$. Let $P_t$ have an invariant probability measure
$\mu$.
\begin{longlist}[(2)]
\item[(1)] If there exists
$r> K^+/ \lambda^2$ such that $\mu({e}^{r|\cdot|^2})<\infty$, then $P_t$
is hypercontractive, that is, $\|P_t\|_{L^2(\mu)\to L^4(\mu)}=1$ holds
for some $t>0$.
\item[(2)] If $\mu({e}^{r|\cdot|^2})<\infty$ holds for all $r>0$, then
$P_t$ is supercontractive,
that is, $\|P_t\|_{L^2(\mu)\to L^4(\mu)}<\infty$
holds for all $t>0$.
\item[(3)] If $P_t {e}^{r|\cdot|^2}$ is bounded for any $t,r>0$, then
$P_t$ is ultracontractive,
that is, $\|P_t\|_{L^2(\mu)\to L^\infty(\mu)}<\infty$
for any $t>0$.
\end{longlist}
\end{cor}
\begin{Remark}\label{Rem11}
To see that results in Corollary~\ref{C1.2} are
sharp, let
$P_t$ be symmetric w.r.t. $\mu$. Then the hypercontractivity is
equivalent to
the validity of the log-Sobolev inequality
\[
\mu(f^2\log f^2)\le C\mu(\GG(f,f)),\qquad f\in C_b^\infty(\R^d),
\mu(f^2)=1,
\]
for some constant $C>0$.
Moreover, if
there exists a constant $R>0$ such that
%
%
\begin{equation}\label{DD} \GG(f,f)\le R^2 |\nn f|^2,\qquad f\in C^\infty
(\R^d),
\end{equation}
we have
$\rr_a(x,y)\ge R^{-1}|x-y|$. So, by the concentration
of measure for the log-Sobolev inequality,
the hypercontractivity implies $\mu({e}^{r|\cdot|^2})<\infty$ for some
$r>0$, while the supercontractivity
implies $\mu({e}^{r|\cdot|^2})<\infty$ for all $r>0$.
Combining this with Corollary~\ref{C1.2}, we have the following
assertions under conditions (A1)--(A3) and (\ref{DD}):
\begin{longlist}[(iii)]
\item[(i)] Let $K\le0$. Then $P_t$ is hypercontractive if and only if
$\mu({e}^{r|\cdot|^2})<\infty$
holds for some $r>0;$
\item[(ii)] $P_t$ is supercontractive if and only if $\mu
({e}^{r|\cdot
|^2})<\infty$ holds for all $r>0;$
\item[(iii)] $P_t$ is ultracontractive if and only if $P_t
{e}^{r|\cdot
|^2}$ is bounded for any $t,r>0$.
\end{longlist}
Therefore, conditions in Corollaries~\ref{C1.2}(2) and~\ref{C1.2}(3)
are sharp for the supercontractivity and ultracontractivity of
$P_t$. Moreover, as shown in~\cite{CW} that when $\si$ is constant,
the sufficient condition $\mu({e}^{r|\cdot|^2})<\infty$ for some $r>
K^+/\lambda^2$ is optimal for the hypercontractivity of $P_t$. So,
Corollary~\ref{C1.2}(1) also provides a sharp sufficient condition
for the hypercontractivity of $P_t$.
\end{Remark}

We will prove Theorem~\ref{T1.1} and Corollary~\ref{C1.2} in the
next section. In Section~\ref{sec3}, we extend these results to SDEs on
Riemannian manifolds possibly with a convex boundary. Finally,
combining results in Section~\ref{sec3} with a conformal change method
introduced in
\cite{W07}, we are able to establish Harnack inequalities in Section
\ref{sec4} for the Neumann semigroup on a class of nonconvex manifolds.

\section{\texorpdfstring{Proofs of Theorem \protect\ref{T1.1} and Corollary \protect
\ref{C1.2}}{Proofs of Theorem 1.1 and Corollary 1.3}}\label{sec2}

Let $x, y\in\R^d, T>0$ and $p> (1+\delta_T/\lambda_T)^2$ be fixed such
that $x\ne y$. We have
%
%
\begin{equation}\label{2.0} \theta_T:= \ff{2\delta_T} {(\sqrt p
-1)\lambda_T} \in
(0,2).
\end{equation}
For $\theta\in(0,2)$, let
\[
\xi_t= \ff{ 2-\theta}{K_T} \bigl(1-{e}^{K_T(t-T)}\bigr),\qquad t\in[0, T].
\]
Then $\xi$ is smooth and strictly positive on $[0, T)$ such that
%
%
\begin{equation}\label{2.1}
2- K_T\xi_t +\xi_t'=\theta,\qquad t\in[0,T].
\end{equation}
Consider the coupling
%
%
\begin{eqnarray}\label{2.2}
\ddd X_t &=&\si(t,X_t)\,\ddd B_t + b(t,X_t)\,\ddd t,\qquad
X_0=x,\nonumber\\
\ddd Y_t &=& \si(t, Y_t)\,\ddd B_t + b(t,Y_t)\,\ddd t \\
&&{}+ \ff1 {\xi_t} \si(t,
Y_t)\si(t,X_t)^{-1} (X_t-Y_t)\,\ddd t,\qquad
Y_0=y.
\nonumber
\end{eqnarray}
Since the additional drift term
$\xi_t^{-1} \si(t, y)\si(t,x)^{-1} (x-y)$ is locally Lipschitzian
in $y$ if (A4) holds, and continuous in $y$ when $\si$ and $b$ are
deterministic and time independent, the coupling $(X_t, Y_t)$ is
a~well-defined continuous process for \mbox{$t<T\land\zeta$}, where $\zeta$ is
the explosion
time of $Y_t$; namely, $\zeta=\lim_{n\to\infty} \zeta_n$ for
\[
\zeta_n:= \inf\{t\in[0,T)\dvtx|Y_t|\ge n\},
\]
where we set $\inf\varnothing=T$.
Let
\[
\ddd\tilde B_t = \ddd B_t +\ff1 {\xi_t} \si(t,X_t)^{-1} (X_t-Y_t)
\,\ddd t,\qquad t<T\land\zeta.
\]
If $\zeta=T$ and
\begin{eqnarray*}
R_s&:=& \exp\biggl[-\int_0^s \xi_t^{-1} \langle\si(t,X_t)^{-1}
(X_t-Y_t),\ddd B_t\rangle\\
&&\hspace*{20pt}{}- \ff1 2 \int_0^s \xi_t^{-2} |\si
(t,X_t)^{-1}(X_t-Y_t)|^2
\,\ddd t\biggr]
\end{eqnarray*}
is a uniformly integrable martingale for $s\in[0,T)$, then by the
martingale convergence theorem, $R_T:=\lim_{t\uparrow T} R_t$
exists and $\{R_t\}_{t\in[0,T]}$ is a martingale. In this case, by
the Girsanov theorem $\{\tilde B_t\}_{t\in[0,T)}$ is a $d$-dimensional
Brownian motion under the probability $R_T \mathbb P$. Rewrite (\ref{2.2})
as
%
%
\begin{eqnarray}\label{2.3}
\ddd X_t &=&\si(t,X_t)\,\ddd\tilde B_t + b(t,X_t)\,\ddd
t-\ff{X_t-Y_t} {\xi_t} \,\ddd t ,\qquad X_0=x,\nonumber\\[-8pt]\\[-8pt]
\ddd Y_t &=& \si(t, Y_t)\,\ddd\tilde B_t + b(t,Y_t)\,\ddd t,
\qquad Y_0=y.
\nonumber
\end{eqnarray}
Since $\int_0^T \xi_t^{-1}\,\ddd
t=\infty$, we will see that the additional drift $-\ff{X_t-Y_t}
{\xi_t} \,\ddd t$ is strong enough to force the coupling to be
successful up
to time $T$. So, we first prove the uniform integrability of
$\{R_{s\land\zeta}\}_{s\in[0,T)}$ w.r.t. $\mathbb P$ so that
$R_{T\land\zeta}:=\lim_{s\uparrow T} R_{s\land\zeta}$ exists, then
prove that $\zeta=T$ $\Q$-a.s.
for $\Q:=R_{T\land\zeta}\mathbb P$ so that $\Q=R_T\mathbb P$.

Let
\[
\tau_n=\inf\{t\in[0,T)\dvtx|X_t|+|Y_t|\ge n\}.
\]
Since $X_t$ is nonexplosive as assumed, we have $\tau_n\uparrow
\zeta$ as $n\uparrow\infty$.
\begin{lem}\label{L2.0} Assume \textup{(A1)}
and \textup{(A2)}. Let $\theta\in(0,2),x,y\in\R^d$ and $T>0$ be fixed.
\begin{longlist}[(2)]
\item[(1)] There holds
\[
\sup_{s\in[0,T), n\ge1} \E R_{s\land\tau_n}\log R_{s\land\tau_n}
\le
\ff{K_T|x-y|^2}{2\lambda_T^2 \theta
(2-\theta)(1-{e}^{-K_TT})|}.
\]
Consequently,
\[
R_{s\land
\zeta}:=\lim_{n\uparrow\infty} R_{s\land\tau_n\land(T-1/n)},\qquad
s\in
[0,T],\qquad
R_{T\land\zeta}:=\lim_{s\uparrow T} R_{s\land\zeta}
\]
exist such
that $\{R_{s\land\zeta}\}_{s\in[0,T]}$ is a uniformly integrable
martingale.
\item[(2)] Let $\Q= R_{T\land\zeta} \mathbb P$. Then
$\Q(\zeta= T)=1$ so that $\Q=R_T\mathbb P$.
\end{longlist}
\end{lem}
\begin{pf} (1) Let $s\in[0,T)$ be fixed. By (\ref{2.3}), (A1)
and the It\^o formula,
\begin{eqnarray*}
\ddd\|X_t-Y_t\|^2 &\le&2\bigl\langle\bigl(\si(t,X_t)-\si(t,Y_t)\bigr)(X_t-Y_t),
\ddd
\tilde
B_t\bigr\rangle\\
&&{}+ K_T|X_t-Y_t|^2\,\ddd t -\ff2 {\xi_t} |X_t-Y_t|^2\,\ddd t
\end{eqnarray*}
holds
for $t\le s\land\tau_n$. Combining this with (\ref{2.1}) we obtain
%
%
\begin{eqnarray}\label{2.5}
\ddd\ff{|X_t-Y_t|^2}{\xi_t}
&\le&\ff2
{\xi_t} \bigl\langle\bigl(\si(t,X_t)-\si(t,Y_t)\bigr)(X_t-Y_t),\ddd\tilde
B_t\bigr\rangle\nonumber\\
&&{}
-\ff{|X_t-Y_t|^2}{\xi_t^2} ( 2-K_T\xi_t +\xi_t')\,\ddd
t\nonumber\\[-8pt]\\[-8pt]
&=&\ff2 {\xi_t} \bigl\langle\bigl(\si(t,X_t)-\si(t,Y_t)\bigr)(X_t-Y_t),\ddd\tilde
B_t\bigr\rangle\nonumber\\
&&{}-\ff
\theta{\xi_t^2} |X_t-Y_t|^2 \,\ddd t,\qquad t\le
s\land\tau_n.
\nonumber
\end{eqnarray}
Multiplying by $\ff1 \theta$ and integrating from $0$ to $s\land\tau_n$,
we obtain
\begin{eqnarray*}
\int_0^{s\land\tau_n} \ff{|X_t-Y_t|^2}{\xi_t^2} \,\ddd t
&\le&\int_0^{s\land\tau_n} \ff2 {\theta\xi_t} \bigl\langle\bigl(\si
(t,X_t)-\si
(t,Y_t)\bigr)(X_t-Y_t),\ddd\tilde B_t\bigr\rangle\\
&&{}-\ff{|X_t-Y_t|^2}{\theta\xi_t} +\ff{|x-y|^2}{\theta\xi_0}.
\end{eqnarray*}
By the Girsanov
theorem, $\{\tilde B_t\}_{t\le\tau_n\land s}$ is the $d$-dimensional
Brownian motion under the probability measure $R_{s\land\tau
_n}\mathbb
P$. So, taking expectation $\E_{s,n}$ with respect to $R_{s\land\tau
_n}\mathbb P$, we arrive at
%
%
\begin{equation}\label{2.00} \E_{s,n} \int_0^{s\land\tau_n}
\ff{|X_t-Y_t|^2}{\xi_t^2}\,\ddd t\le\ff{|x-y|^2}{\theta\xi_0},\qquad s\in
[0,T), n\ge1.
\end{equation}
By (A2) and the definitions of $R_t$
and $\tilde B_t$, we have
\begin{eqnarray}
\log R_r &=& -\int_0^r\ff1 {\xi_t}
\langle\si(t,X_t)^{-1} (X_t-Y_t), \ddd\tilde B_t\rangle+\ff1 2
\int
_0^r \ff
{|\si(t,X_t)^{-1} (X_t,Y_t)|^2}{\xi_t^2}\,\ddd t\nonumber\\
&\le&-\int_0^r\ff1 {\xi_t} \langle\si(t,X_t)^{-1} (X_t-Y_t), \ddd
\tilde
B_t\rangle+\ff1 {2\lambda_T^2} \int_0^r\ff{|X_t-Y_t|^2}{\xi
_t^2}\,\ddd
t,\nonumber\\
&&\eqntext{r\le s\land\tau_n.}
\end{eqnarray}
Since $\{\tilde B_t\}$ is
the $d$-dimensional Brownian motion under $R_{s\land\tau_n}\mathbb P$ up
to $s\land\tau_n$, combining this with (\ref{2.00}), we obtain
\[
\E R_{s\land\tau_n}\log R_{s\land\tau_n} =\E_{s,n} \log R_{s\land
\tau_n} \le\ff{|x-y|^2}{2\lambda_T^2\theta\xi_0},\qquad s\in[0,T),
n\ge
1.
\]
By the martingale convergence theorem and the Fatou lemma,
$\{R_{s\land\zeta}\dvtx s\in
[0,T]\}$ is a well-defined martingale with
\[
\E R_{s\land\zeta}\log R_{s\land\zeta} \le\ff{|x-y|^2}{2\lambda
_T^2\theta\xi_0}=
\ff{K_T|x-y|^2}{2\lambda_T^2 \theta(2-\theta) (1-{e}^{-K_TT})},\qquad
s\in
[0,T].
\]
To see that $\{R_{s\land\zeta}\dvtx s\in[0,T]\}$ is a martingale, let
$0\le s<t\le T$. By the dominated convergence theorem and the
martingale property of
$\{R_{s\land\tau_n}\dvtx s\in[0,T)\}$, we have
\begin{eqnarray*}
\E(R_{t\land\zeta}|\F_s)&=&\E\Bigl(\lim_{n\to\infty
}R_{t\land\tau_n\land(T-1/n)}|\F_s\Bigr)
= \lim_{n\to\infty} \E\bigl(R_{t\land\tau_n\land(T-1/n)}|\F_s\bigr) \\
&=&\lim_{n\to\infty} R_{s\land\tau_n} = R_{s\land\zeta}.
\end{eqnarray*}

(2) Let $\si_n=\inf\{t\ge0\dvtx|X_t|\ge n\}$. We have
$\si_n\uparrow\infty$ $\mathbb P$-a.s and hence, also $\Q$-a.s. Since
$\{\tilde B_t\}$ is a $\Q$-Brownian motion up to $T\land\zeta$, it
follows from (\ref{2.5}) that
\[
\ff{(n-m)^2}{\xi_0} \Q(\si_m>t, \zeta_n\le t) \le\E_\Q
\ff{|X_{t\land\si_m\land\zeta_n} - X_{t\land\si_m\land
\zeta_n}|^2}{\xi_{t\land\si_m\land\zeta_n}} \le
\ff{|x-y|^2}{\xi_0}
\]
holds for all $n>m>0$ and $t\in[0,T)$. By
letting first $n\uparrow\infty$ then $m\uparrow\infty$, we obtain
$\Q(\zeta\le t)=0$ for all $t\in[0,T)$. This is equivalent to
$\Q(\zeta=T)=1$ according to the definition of $\zeta$.
\end{pf}

Lemma~\ref{L2.0} ensures that under $\Q:= R_{T\land\zeta}\mathbb P$,
$\{\tilde B_t\}_{t\in[0,T]}$ is a Brownian motion. Then by (\ref
{2.3}), the
coupling $(X_t,Y_t)$ is well-constructed under $\Q$ for $t\in[0,T]$.
Since $\int_0^T\xi_t^{-1}\,\ddd t=\infty$, we shall see
that the coupling is successful up to time $T$, so that $X_T=Y_T$
holds $\Q$-a.s. (see the proof of Theorem~\ref{T1.1} below). This
will provide
the desired Harnack inequality for $P_t$ as explained in Section~\ref{sec1}
as soon as $R_{T\land\zeta}$ has finite $p/(p-1)$-moment. The next
lemma provides an explicit upper bound on moments of
$R_{T\land\zeta}$.
\begin{lem}\label{L2.1} Assume \textup{(A1)--(A3)}. Let $R_t$ and $\xi
_t$ be
fixed for $\theta=\theta_T$. We have
%
%
\begin{eqnarray}\label{2.4}
&&\sup_{s\in[0,T]} \E\biggl\{R_{s\land\zeta}
\exp\biggl[\ff{\theta_T^2}{8\delta_T^2} \int_0^{s\land\zeta}
\ff{|X_t-Y_t|^2}{\xi_t^2}\,\ddd t\biggr]\biggr\}\nonumber\\[-8pt]\\[-8pt]
&&\qquad\le\exp\biggl[ \ff{\theta_T
K_T |x-y|^2}{4\delta_T^2 (2-\theta_T)
(1-{e}^{-K_TT})}\biggr].\nonumber
\end{eqnarray}
Consequently,
%
%
\begin{equation}\label{2.4'} \sup_{s\in[0,T]} \E R_{s\land\zeta
}^{1+r_T} \le
\exp\biggl[\ff{\theta_T K_T(2\delta_T +\theta_T\lambda
_T)|x-y|^2}{8\delta_T^2
(2-\theta_T)(\delta_T +\theta_T \lambda_T)
(1-{e}^{-K_TT})}\biggr]
\end{equation}
holds for
\[
r_T= \ff{\lambda_T^2 \theta_T^2}{4\delta_T^2 + 4 \theta_T \lambda
_T \delta_T}.
\]
\end{lem}
\begin{pf} Let $\theta=\theta_T$. By (\ref{2.5}), for any $r>0$ we
have
\begin{eqnarray*}
&&\E_{s,n} \exp\biggl[r \int_0^{s\land\tau_n} \ff
{|X_t-Y_t|^2}{\xi_t^2}\,\ddd t\biggr]\\
&&\qquad\le\exp\biggl[\ff{r|x-y|^2}{\theta_T\xi_0}\biggr]\\
&&\qquad\quad{}\times\E_{s,n}
\exp\biggl[\ff{2r}{\theta_T} \int_0^{s\land\tau_n}
\ff1 {\xi_t} \bigl\langle\bigl(\si(t,X_t)-\si(t,Y_t)\bigr)(X_t-Y_t), \ddd\tilde
B_t\bigr\rangle\biggr]\\
&&\qquad\le\exp\biggl[\ff{rK_T|x-y|^2}{\theta_T(2-\theta_T)
(1-{e}^{-K_TT})}\biggr]\\
&&\qquad\quad{}\times\biggl(\E_{s,n}
\exp\biggl[\ff{8r^2\delta_T^2}{\theta_T^2}\int_0^{s\land\tau_n}
\ff{|X_t-Y_t|^2}{\xi_t^2}\,\ddd
t\biggr]\biggr)^{1/2},
\end{eqnarray*}
where the last step
is due to (A3) and the fact that
\[
\E{e}^{M_t}\le\bigl(\E{e}^{2\langle M\rangle_t}\bigr)^{1/2}
\]
for a continuous exponential integrable martingale $M_t$.
Taking $r= \theta_T^2/(8\delta_T^2)$, we arrive at
\[
\E_{s,n} \exp\biggl[\ff{\theta_T^2 }{8\delta_T^2}\int_0^{s\land\tau_n}
\ff{|X_t-Y_t|^2}{\xi_t^2}\,\ddd t\biggr]\le\biggl[\ff{\theta_T
K_T|x-y|^2}{4\delta_T^2(2-\theta_T) (1-{e}^{-K_TT})}\biggr],\qquad n\ge1.
\]
This implies (\ref{2.4}) by letting $n\to\infty$.

Next, by (A2) and the definition of $R_s$, we have
%
%
\begin{eqnarray}\label{W0}
\E R_{s\land\tau_n}^{1+r_T}
&=& \E_{s,n} R_{s\land\tau_n}^{r_T}\nonumber\\
&=&
\E_{s,n}\exp\biggl[- r_T\int_0^{s\land\tau_n} \ff{1}{\xi_t}
\langle\si(t,X_t)^{-1}(X_t-Y_t),\ddd\tilde B_t\rangle\\
&&\hspace*{51.7pt}{} + \ff{r_T}2 \int_0^{s\land\tau_n} \ff{
|\si(t,X_t)^{-1} (X_t-Y_t)|^2}{\xi_t^2}
\,\ddd t\biggr].
\nonumber
\end{eqnarray}
Noting that for any exponential integrable martingale $M_t$
w.r.t. $R_{s\land\tau_n}\mathbb P$, one has
\begin{eqnarray*}
&&\E_{s,n}\exp[r_TM_t +r_T\langle M\rangle_t/2]\\[-2pt]
&&\qquad= \E_{s,n}
\exp[r_TM_t -r_T^2q\langle M\rangle_t/2 + r_T(qr_T+1)\langle M\rangle
_t/2]\\[-2pt]
&&\qquad\le(\E_{s,n} \exp[r_T q M_t- r_T^2q^2\langle M\rangle_t/2]
)^{1/q} \\[-2pt]
&&\qquad\quad\hspace*{0pt}{}\times\biggl(\E_{s,n} \exp\biggl[\ff{r_Tq (r_Tq+1)}{2(q-1)} \langle
M\rangle_t\biggr]\biggr)^{(q-1)/q}\\[-2pt]
&&\qquad= \biggl(\E_{s,n} \exp\biggl[\ff{r_Tq (r_Tq+1)}{2(q-1)}
\langle M\rangle_t\biggr]\biggr)^{(q-1)/q},\qquad q>1,
\end{eqnarray*}
it follows from (\ref{W0}) that
%
%
\begin{equation}\label{2.6}\qquad
\E R_{s\land\tau_n}^{1+r_T} \le\biggl(\E
_{s,n} \exp
\biggl[ \ff{qr_T(qr_T+1)}{2(q-1)\lambda_T^2}
\int_0^{s\land\tau_n} \ff{| X_t-Y_t|^2}{\xi_t^2} \,\ddd
t\biggr]\biggr)^{(q-1)/q}.
\end{equation}
Take
%
%
\begin{equation}\label{2.01}q= 1 +\sqrt{1+ r_T^{-1}},
\end{equation}
which minimizes
$q(qr_T+1)/(q-1)$ such that
%
%
\begin{eqnarray}\label{2.02}
\ff{qr_T(qr_T+1)}{2\lambda_T^2(q-1)}&=& \ff{r_T+\sqrt
{r_T(r_T+1)}}{2\lambda_T^2
\sqrt{1+r_T^{-1}}}\bigl(r_T+1+\sqrt{r_T(r_T+1)}\bigr)\nonumber\\[-9pt]\\[-9pt]
&=& \ff{(r_T+ \sqrt{r_T^2 +r_T})^2}{2\lambda_T^2}=
\ff{\theta_T^2}{8\delta_T^2}.
\nonumber
\end{eqnarray}
Combining
(\ref{2.6}) with (\ref{2.4}) and (\ref{2.02}), and noting that due
to (\ref{2.01}) and the definition of $r_T$
\[
\ff{q-1} q = \ff{\sqrt{1+r_T^{-1}}}{1+ \sqrt{1+r_T^{-1}}} = \ff
{2\delta_T
+\theta_T\lambda_T}{2\delta_T +2\theta_T\lambda_T},
\]
we obtain
\[
\E R_{s\land\tau_n}^{1+r_T} \le\exp\biggl[\ff{\theta_T K_T(2\delta_T
+\theta
_T\lambda_T)|x-y|^2}
{8\delta_T^2 (2-\theta_T)(\delta_T +\theta_T \lambda_T)
(1-{e}^{-K_TT})}\biggr].
\]
According to the Fatou lemma, the proof is
then completed by letting $n\to\infty$.
\end{pf}
\begin{pf*}{Proof of Theorem~\ref{T1.1}} Since (A3) also
holds for $\delta_{p,T}$ in place of $\delta_T$, it suffices to prove the
desired Harnack inequality for $\delta_T$ in place of $\delta_{p,T}$.

(1) By Lemma~\ref{L2.0}, $\{R_{s\land\zeta}\}_{s\in[0,T]}$ is an
uniformly integrable martingale and $\{\tilde B_t\}_{t\le T}$ is a
$d$-dimensional Brownian motion\vadjust{\goodbreak} under the probability $\Q$. Thus,
$Y_t$ can be solved up to time $T$. Let
\[
\tau=\inf\{t\in[0,T]\dvtx X_t=Y_t\}
\]
and set $\inf\varnothing=\infty$ by convention. We claim that $\tau
\le
T$ and thus, $X_T=Y_T$, $\Q$-a.s.
Indeed, if for some $\omega\in\OO$ such that $\tau(\omega)>T$, by
the continuity of the processes we have
\[
\inf_{t\in[0,T]} |X_t-Y_t|^2(\omega)>0.
\]
So,
\[
\int_0^T \ff{|X_t-Y_t|^2}{\xi_t^2}\,\ddd t=\infty
\]
holds on the set $\{\tau>T\}$.
But according to Lemma~\ref{L2.1}, we have
\[
\E_\Q\int_0^T\ff{|X_t-Y_t|^2}{\xi_t^2}\,\ddd
t <\infty,
\]
we conclude that $\Q(\tau>T)=0$. Therefore, $X_T=Y_T$
$\Q$-a.s.

Now, combining Lemma~\ref{L2.0} with $X_T=Y_T$ and using the Young
inequality, for $f\ge1$ we have
\begin{eqnarray*}
P_T\log f(y)&=& \E_\Q[\log f(Y_T)]= \E
[R_{T\land\zeta}\log f(X_T)] \\
&\le&\E R_{T\land\zeta}\log R_{T\land\zeta} +\log\R f(X_T)\\
&\le&\log
P_Tf(x)+ \ff{K_T|x-y|^2}{2\lambda_T^2\theta
(2-\theta)(1-{e}^{-K_TT})}.
\end{eqnarray*}
This completes
the proof of (1) by taking $\theta=1$.

(2) Let $\theta=\theta_T$. Since $X_T=Y_T$ and $\{\tilde B_t\}_{t\in
[0,T]}$ is the $d$-dimensional Brownian motion under $\Q$, we have
%
%
\begin{eqnarray}\label{2.8}
(P_Tf(y))^p&=& (\E_\Q[f(Y_T)])^p =(\E
[R_{T\land\zeta}
f(X_T)])^p\nonumber\\[-8pt]\\[-8pt]
&\le&(P_Tf^p(x))\bigl(\E
R_{T\land\zeta}^{p/(p-1)}\bigr)^{p-1}.\nonumber
\end{eqnarray}
Due to (\ref{2.0}),
we see that
\[
\ff p{p-1} = 1 + \ff{\lambda_T^2 \theta_T^2}{4\delta_T(\delta_T+
\theta
_T\lambda_T)}.
\]
So,
it follows from Lemma~\ref{L2.1} and (\ref{2.0}) that
\begin{eqnarray*}
\bigl(\E R_{T\land\zeta}^{p/(p-1)}\bigr)^{p-1} &=&
(\E R_{T\land\zeta}^{1+r_T})^{p-1}\le\exp\biggl[\ff{(p-1)\theta_T
K_T(2\delta_T+\theta_T\lambda_T)|x-y|^2}{8\delta_T^2(2-\theta
_T)(\delta_T+\theta
_T\lambda_T)(1-{e}^{-K_TT})}\biggr]\\
&=& \exp\biggl[\ff{K_T\sqrt p (\sqrt p-1)|x-y|^2}{4\delta_T[(\sqrt
p-1)\lambda_T-\delta_T] (1-{e}^{-K_TT})}\biggr].
\end{eqnarray*}
Then the
proof is finished by combining this with (\ref{2.8}).\vadjust{\goodbreak}
\end{pf*}
\begin{pf*}{Proof of Corollary~\ref{C1.2}} Let $f\in
\B_b^+(\R^d)$ be such that $\mu(f^p)\le1$. Let $p>(1+\delta
/\lambda)^2$.
By Theorem~\ref{T1.1}(2), we have
\[
(P_tf(y))^p\exp\biggl[-\ff{K\sqrt p (\sqrt p-1)|x-y|^2}{4\delta_p [(\sqrt
p-1)\lambda-\delta_p] (1-{e}^{-Kt})}\biggr]
\le P_t f^p(x), \qquad x,y\in\R^d,
\]
where $\delta_p=\max\{\delta, \ff\lambda2
(\sqrt p-1)\}$. Integrating w.r.t. $\mu(\ddd x)$ and noting that $\mu
$ is
$P_t$-invariant, we obtain
%
%
\begin{equation}\label{2.9}\quad
(P_tf(y))^p \int_{\R^d} \exp\biggl[-\ff {K\sqrt p (\sqrt
p-1)|x-y|^2}{4\delta[(\sqrt p-1)\lambda-\delta]
(1-{e}^{-Kt})}\biggr]\mu (\ddd x) \le1.
\end{equation}
Taking $f= n\land(p_t(y,\cdot))^{1/p}$ and
letting $n\uparrow\infty$, we prove the first assertion.

Next, let $B(0,1)=\{x\in\R^d\dvtx|x|\le1\}$. Since $\mu$ is an
invariant measure, it has a strictly positive density w.r.t. the
Lebesgue measure so that $\mu(B(0,1))>0$ (cf.~\cite{BKR}). Let $p\ge
(1+ 2\delta/\lambda)^2$. We have $\delta_p= (\sqrt p-1)\lambda/2$
and thus
\[
\ff{ \sqrt p (\sqrt p-1) }{4\delta_p [(\sqrt p-1)\lambda-\delta_p]
}= \ff
{\sqrt p}{\lambda^2(\sqrt p-1)}.
\]
Combining this with (\ref{2.9}) and noting that
\begin{eqnarray*}
&&\int_{\R^d} \exp\biggl[-\ff{K\sqrt p (\sqrt
p-1)|x-y|^2}{4\delta[(\sqrt p-1)\lambda-\delta]
(1-{e}^{-Kt})}\biggr]\mu(\ddd x) \\
&&\qquad\ge\mu(B(0,1))\exp\biggl[-\ff{K\sqrt p (\sqrt p-1)(1+|y|)^2}{4\delta
[(\sqrt p-1)\lambda-\delta] (1-{e}^{-Kt})}\biggr],
\end{eqnarray*}
we obtain
%
%
\begin{equation}\label{2.10}\quad
(P_t f(y))^p \le C_1 \exp\biggl[\ff{K\sqrt p
(1+|y|)^2}{\lambda_T^2(\sqrt p-1) (1-{e}^{-Kt})}\biggr],\qquad t>0, y\in
\R^d,
\end{equation}
for some constant $C_1>0$ and all $f\in
\B_b^+(\R^d)$ with $\mu(f^p)\le1$. Since
\[
\lim_{p\to\infty} \lim_{t\to\infty} \ff{K\sqrt p }{\lambda
^2(\sqrt
p-1)(1-{e}^{-Kt})}=\ff{K^+}{\lambda^2},
\]
for any
$r> K^+/\lambda^2$ there exist $p>(1+2\delta_T/\lambda)^2, \bb>1$
and $t_1>0$
such that
\[
(P_{t_1}f(y))^{\bb p} \le C_2 {e}^{r|y|^2},\qquad y\in\R^d, f\in\B
_b^+(\R^d), \mu(f^p)\le1,
\]
holds for some constant $C_2>0$. Thus, $\mu({e}^{r|\cdot|^2})
<\infty$ implies that
\[
\|P_{t_1}\|_{L^p(\mu)\to
L^{p\bb}(\mu)}<\infty.
\]
Since $\|P_s\|_{L^q(\mu)}=1$ holds for any
$q\in[1,\infty]$, by the interpolation theorem and the semigroup
property one may find $t_2>t_1$ such that
%
%
\begin{equation}\label{2.11} \|P_{t_2}\|_{L^2(\mu)\to L^4(\mu
)}<\infty.
\end{equation}
Moreover, by~\cite{H}, Theorem 3.6(ii), there exist some constants
$\eta
, C_3>0$
such that
\[
\|P_t-\mu\|_{L^2(\mu)}\le C_3 {e}^{-\eta t},\qquad t\ge0.
\]
Combining this with (\ref{2.11}) we conclude that $\|P_t\|_{L^2(\mu
)\to
L^4(\mu)}\le1$
holds for sufficiently large $t>0$, that is, (2) holds.

Finally, (3) and (4) follow immediately from (\ref{2.10}) and the
interpolation theorem.
\end{pf*}

\section{Extension to manifolds with convex boundary}\label{sec3}

Let $M$ be a $d$-dimen\-sional complete, connected Riemannian
manifold, possibly with a convex boundary $\pp M$. Let $N$ be the
inward unit normal vector filed of $\pp M$ when $\pp M\ne\varnothing$.
Let $P_t$ be the (Neumann) semigroup generated by
\[
L:= \psi^2(\DD+Z)
\]
on $M$, where $\psi\in C^1(M)$ and $Z$ is a $C^1$ vector field on
$M$. Assume that $\psi$ is bounded and
%
%
\begin{equation}
\Ric-\nn Z\ge-K_0
\end{equation}
holds for some
constant $K_0\ge0$. Then the (reflecting) diffusion process
generated by $L$ is nonexplosive.

To formulate $P_t$ as the semigroup associated to a SDE like
(\ref{1.1}), we set
%
%
\begin{equation}
\label{3.0} \si= \sqrt2 \psi,\qquad b= \psi^2 Z.
\end{equation}
Let $\ddd_I$ denote the It\^o differential on $M$. In local coordinates
the It\^o differential
for a continuous semi-martingale $X_t$ on $M$ is given by (see~\cite{ATW}
or~\cite{E})
\[
(\ddd_I X_t)^k =\ddd X_t^k +\ff1 2 \sum_{i,j=1}^d \GG
_{ij}^k(X_t)\,\ddd
\langle X^i,X^j\rangle_t,\qquad 1\le k\le d.
\]
Then $P_t$ is the semigroup for the solution to the
SDE
%
%
\begin{equation}\label{3.1} \ddd_I X_t= \si(X_t)\Phi_t\,\ddd B_t +
b(X_t) \,\ddd t
+N(X_t)\,\ddd l_t,
\end{equation}
where $B_t$ is the $d$-dimensional
Brownian motion on a complete filtered probability space $(\OO,
\{\F_t\}_{t\ge0}, \mathbb P)$, $\Phi_t$ is the horizontal lift of $X_t$
onto the frame bundle $O(M)$, and $l_t$ is the local time of $X_t$
on $\pp M$. When $\pp M=\varnothing$, we simply set $l_t=0$.

To derive the Harnack inequality as in Section~\ref{sec2}, we assume that
%
%
\begin{equation}\label{3.2} \lambda:=\inf\si>0, \qquad \delta:=
\sup\si-\inf\si<\infty.
\end{equation}
Now, let $x,y\in M$ and $T>0$
be fixed. Let $\rr$ be the Riemannian distance on $M$, that is,
$\rr(x,y)$ is the length of the minimal geodesic on $M$ linking $x$
and $y$, which exits if $\pp M$ is either convex or empty.

Let
$X_t$ solve (\ref{3.1}) with $X_0=x$. Next, any strictly positive
function $\xi\in C([0,T))$, let $Y_t$ solve
\begin{eqnarray*}
\ddd_I Y_t&=& \si(Y_t)P_{X_t,Y_t}\Phi_t\,\ddd B_t + b(X_t) \,\ddd t\\
&&{}-\ff
{\si
(Y_t)\rr(X_t,Y_t)}{\si(X_t)\xi_t}
\nn\rr(X_t, \cdot)(Y_t)\,\ddd t +N(Y_t)\,\ddd\tilde l_t
\end{eqnarray*}
for $Y_0=y$, where
$\tilde l_t$ is the local time of $Y_t$ on $\pp M$, and
$P_{X_t,Y_t}\dvtx
T_{X_t}M\to T_{Y_t}M$ is the parallel displacement along the minimal
geodesic from $X_t$ to $Y_t$, which exists since $\pp M$ is convex
or empty. As explained in~\cite{ATW}, Section 3, we may and do
assume that the cut-locus of $M$ is empty such that the parallel
displacement is smooth. Let
\[
\ddd\tilde B_t =\ddd B_t +\ff{\rr(X_t,Y_t)}{\xi_t \si(X_t)} \Phi_t^{-1}
\nn\rr(\cdot, Y_t)(X_t)\,\ddd t,\qquad t<T.
\]
By the Girsanov theorem,
for any $s\in(0,T)$ the process $\{\tilde B_t\}_{t\in[0,s]}$ is the
$d$-dimensional Brownian motion under the weighted probability
measure $R_s\mathbb P$, where
%
%
\begin{eqnarray}\label{R}
R_s&:=&
\exp\biggl[-\int_0^s\ff{\rr(X_t,Y_t)}{\xi_t\si(X_t)} \langle\nn
\rr(\cdot,Y_t)(X_t), \Phi_t \,\ddd B_t\rangle\nonumber\\[-8pt]\\[-8pt]
&&\hspace*{98.7pt}{}-\ff1 2 \int_0^s
\ff{\rr(X_t,Y_t)^2}{\xi_t^2 \si(X_t)^2}\,\ddd t\biggr].\nonumber
\end{eqnarray}
Thus, by (\ref{3.0}) we have
\begin{eqnarray*}
\ddd_I X_t&=& \sqrt2 \psi(X_t)\Phi_t\,\ddd\tilde
B_t + (\psi^2Z)(X_t) \,\ddd t\\
&&{}-\ff{\rr(X_t,Y_t)}{\xi_t} \nn\rr(\cdot,
Y_t)(X_t)\,\ddd t
+N(X_t)\,\ddd l_t,\\
\ddd_I Y_t&=& \sqrt2 \psi(Y_t)\Phi_t\,\ddd\tilde B_t + (\psi^2Z)(Y_t)
\,\ddd t
+N(Y_t)\,\ddd\tilde l_t.
\end{eqnarray*}
Let $\xi\in C^1([0,T))$
be strictly positive and take
\[
\bb_t= -\ff{\rr(X_t,Y_t)}{\sqrt2 \xi_t \psi(X_t)}\Phi_t^{-1}
\nn\rr(\cdot, Y_t)(X_t).
\]
Repeating the proof of (4.10) in~\cite{W09},
we obtain
\begin{eqnarray*}
\ddd\rr(X_t,Y_t)&\le&\bigl(\si(X_t)-\si(Y_t)\bigr) \langle\nn\rr(\cdot, Y_t)
(X_t), \Phi_t\,\ddd\tilde B_t\rangle\\
&&{}+ K_1 \rr(X_t,Y_t)\,\ddd t
-\ff{\rr(X_t,Y_t)}{\xi_t}\,\ddd t,\qquad t<T,
\end{eqnarray*}
where
\[
K_1= K_0\|\psi\|_\infty^2 + 2\|Z\|_\infty\|\nn\psi\|_\infty\|\psi
\|
_\infty.
\]
This implies that
\begin{eqnarray*}
\ddd\ff{\rr(X_t, Y_t)^2}{\xi_t}&\le&\ff2 {\xi_t} \rr(X_t,Y_t)
\bigl(\si(X_t)-\si(Y_t)\bigr)\langle\nn\rr(\cdot, Y_t)(X_t), \Phi_t\,\ddd
\tilde
B_t\rangle\\
&&{}-\ff{\rr(X_t,Y_t)^2}{\xi_t^2}(2- K\xi_t + \xi_t')\,\ddd t
\end{eqnarray*}
holds for $t<T$ and
%
%
\begin{eqnarray}\label{K}
K:\!&=& 2K_1 + \|\nn\si\|_\infty^2\nonumber\\[-8pt]\\[-8pt]
&=&2K_0\|\psi\|_\infty^2 + 4\|Z\|_\infty
\|\nn\psi\|_\infty\|\psi\|_\infty+ 2 \|\nn
\psi\|_\infty^2.
\nonumber
\end{eqnarray}
In particular, letting
\[
\xi_t=
\ff{2-\theta}K \bigl(1-{e}^{K(t-T)}\bigr),\qquad t\in[0,T], \theta\in(0,2),
\]
we
have
\[
2-K\xi_t +\xi_t'= \theta.
\]
Therefore, the following result follows
immediately by repeating calculations in Section~\ref{sec2}.
\begin{theorem}\label{T3.1} Assume that $\pp M$ is either empty or convex.
Let (\ref{C}) and $Z,\phi$ be bounded such that
\[
K:=2K_0\|\psi\|_\infty^2 + 4\|Z\|_\infty
\|\nn\psi\|_\infty\|\psi\|_\infty+ 2 \|\nn\psi\|_\infty
^2<\infty.
\]
Then all assertions in Theorem~\ref{T1.1} and Corollaries
\ref{C1.0},~\ref{C1.2} hold for $P_t$ the (Neumann) semigroup
generated by $L= \psi^2(\DD+Z)$ on $M$ with $\rho(x, y)$ replacing $|x-y|$, and for
constant functions $K, \delta:= \sup\psi-\inf\psi$ and $\lambda:=
{\inf}|\psi|$.
\end{theorem}

\section{Neumann semigroup on nonconvex manifolds}\label{sec4}

Following the line of~\cite{W07b}, we are able to make the boundary from
nonconvex to convex by using a conformal
change of metric.
This will enable us to extend our results to the Neumann semigroup on a
class of nonconvex manifolds.

Let $\pp M \ne\varnothing$ with $N$ the inward normal unit vector field.
Then the second fundamental form of $\pp M$
is a two-tensor on the tangent space of $\pp M$ defined by
\[
\II(X,Y):=-\langle\nn_X N, Y\rangle,\qquad X,Y\in T\,\pp M.
\]
Assume that there exists $\kk>0$ and $K_0\in\R$ such that
%
%
\begin{equation}\label{C} \Ric-\nn Z\ge-K_0, \qquad\II\ge-\kk
\end{equation}
holds for $M$ and a $C^1$ vector field $Z$.
We shall consider the Harnack inequality for
the Neumann semigroup $P_t$ generated by
\[
L=\DD+Z.\vadjust{\goodbreak}
\]

To make the boundary convex, let $f\in C_b^\infty(M)$ such that $f\ge
1$ and $N\log f|_{\pp M}\ge\kk$.
By~\cite{W07b}, Lemma 2.1,
$\pp M$ is convex under the metric
\[
\langle\cdot,\cdot\rangle'= f^{-2} \langle\cdot,\cdot\rangle.
\]
Let $\DD'$ and $\nn'$ be the Laplacian and gradient induced
by the new metric.
We have (see (2.2) in~\cite{TW})
\[
L= f^{-2} (\DD'+Z'),\qquad Z'= f^2 Z +\ff{d-2} 2\nn f^2.
\]
Let $\Ric'$ be the Ricci curvature induced by the metric $\langle
\cdot
,\cdot\rangle'$. We have (see the proof of
\cite{W09}, Theorem 5.1)
\[
\Ric'-\nn' Z'\ge-K_f\langle\cdot,\cdot\rangle'
\]
for
%
%
\begin{equation}\label{4.1} K_f=\sup\{Kf^2 -d\DD f +(d-3) |\nn f|^2 + 3
|Z|f|\nn f|\}.
\end{equation}
Applying
Theorem~\ref{T3.1} to the convex manifold $(M,\langle\cdot,\cdot
\rangle')$, $\psi= f^{-1}$ and
%
%
\begin{eqnarray} \label{4.2}
K &=& 2 K_f^+ \|f^{-1}\|_\infty+ 4\|Z'\|'_\infty\|\nn'
f^{-1}\|'_\infty
\|f^{-1}\|_\infty+ 2\|\nn' f^{-1}\|_\infty'^2\nonumber\\[-8pt]\\[-8pt]
&\le&2K_f^++ 4\|fZ+ (d-2)\nn f\|_\infty\|\nn f\|_\infty+ 2\|\nn f\|
_\infty^2,
\nonumber
\end{eqnarray}
where $\|\cdot\|'$ is the norm induced by $\langle\cdot,\cdot
\rangle'$
and we have used that $f\ge1$, we obtain the following result.
\begin{theorem} \label{T4.1}Let (\ref{C}) hold for some $\kk>0$ and
$K_0\in\R$, and let $P_t$ be the Neumann semigroup generated by
$L=\DD+Z$ on $M$. Then for any $f\in C_b^\infty(M)$ such that $\inf
f= 1$, $N\log f|_{\pp M}\ge\kk$ and $K<\infty$, where $K$ is fixed
by (\ref{4.1}) and (\ref{4.2}), all assertions in Theorem
\ref{T1.1} and Corollaries~\ref{C1.0} and~\ref{C1.2} hold
with $\rho(x, y)$ replacing $|x-y|$ for
constant functions $K, \delta:= \sup f^{-1}-\inf f^{-1}$ and $\lambda:=
\inf f^{-1}$.
\end{theorem}
\begin{Remark}\label{Rem41}
A simple choice of $f$ in Theorem~\ref{T4.1}
is $f=\phi\circ\rr_\pp$, where $\rr_\pp$ is the Riemannian
distance to
the boundary which is smooth on $\{\rr_\pp\le r_T\}$ for some
$r_T>0$ provided the injectivity radius of the boundary is positive,
and $f\in C_b^\infty([0,\infty))$ is such that $f(0)=1, f'(0)=\kk$
and $f(r)=f(r_T)$ for $r\ge r_T$. In general, $f$ is taken according to
$r_T$ and bounds of the second fundamental form and sectional
curvatures, see, for example,~\cite{W07b,W09} for details. With specific
choices of~$f$, Theorem~\ref{T4.1} provides explicit Harnack type
inequalities, heat kernels estimates and criteria on contractivity
properties for the Neumann semigroup on manifolds with nonconvex
boundary.
\end{Remark}
%









\section*{Acknowledgment}
The author would like to thank the referees and\break Dr.~Wei~Liu for useful comments
and corrections.

%

%
\printaddresses

\end{document}